\documentclass[12pt]{amsart}
\usepackage{extpfeil}
\usepackage[all]{xy}
\usepackage{amssymb} 
\usepackage[margin=1in]{geometry}
\usepackage{euscript}
\usepackage{mathtools}
\usepackage{amsthm}
\usepackage{mathrsfs}
\usepackage{ytableau}
\usepackage{color}
\usepackage[bookmarks, colorlinks=true, linkcolor=blue, citecolor=blue, urlcolor=blue]{hyperref}

\numberwithin{equation}{section}

\theoremstyle{plain}
	
	\newtheorem{prop}[equation]{Proposition}

\theoremstyle{definition}
	\newtheorem{defn}[equation]{Definition}
	\newtheorem{ex}[equation]{Example}

\theoremstyle{remark}
	\newtheorem{rem}[equation]{Remark}
 
\def\nc{\newcommand}
\def\on{\operatorname}

\nc{\C}{\mathbb{C}}
\nc{\Z}{\mathbb{Z}}
\nc{\PP}{\mathbb{P}}
\nc{\R}{\mathbb{R}}

\nc{\id}{{\on{id}}}
\nc\Hom{{\on{Hom}}}
\nc\cone{{\on{cone}}}
\nc\Ob{{\on{Ob}}}
\nc\Spec{{\on{Spec}}}
\nc\Mod{{\on{Mod}}}
\nc\Perf{{\on{Perf}}}
\nc\End{{\on{End}}}
\nc{\into}{\hookrightarrow}
\nc{\tr}{\on{tr}}
\nc{\im}{\on{im}}
\nc{\pt}{\on{pt}}
\nc{\coker}{\on{coker}}
\nc{\rk}{\on{rank}}
\nc{\gr}{\on{gr}}
\nc{\Sym}{\on{Sym}}
\nc{\xra}{\xrightarrow}
\nc{\Feven}{F_{\on{even}}}
\nc{\Fodd}{F_{\on{odd}}}
\nc{\la}{\leftarrow}
\nc{\xla}{\xleftarrow}
\nc{\onto}{\twoheadrightarrow}

\def \l {\lambda}
\def \S {\on{ \bf S}}
\def \th {\on{th}}

\title{Computing Schur Complexes}

\author{Michael K. Brown}
\author{Hang Huang}
\author{Robert P. Laudone}
\author{Michael Perlman}
\author{Claudiu Raicu}
\author{Steven V Sam}
\author{Jo\~{a}o Pedro Santos}

\newcommand{\Addresses}{{
	\vskip\baselineskip
  	\footnotesize
  	\noindent \textsc{Department of Mathematics, University of Wisconsin-Madison} \par\nopagebreak
	\noindent \textit{E-mail addresses:} \texttt{mkbrown5@wisc.edu, hhuang235@math.wisc.edu, laudone@wisc.edu}
  	\vskip\baselineskip
  	\noindent \textsc{Department of Mathematics, University of Notre Dame} \par\nopagebreak
	\noindent \textit{E-mail addresses:} \texttt{mperlman@nd.edu, craicu@nd.edu, jsantos3@nd.edu}
  	\vskip\baselineskip
  	\noindent \textsc{Department of Mathematics, University of California, San Diego} \par\nopagebreak
	\noindent \textit{E-mail addresses:} \texttt{ssam@ucsd.edu}
}}

\renewcommand{\thefootnote}

\begin{document}
\pagestyle{plain}
\maketitle
\footnote{The authors gratefully acknowledge support from the National Science Foundation (NSF award DMS-1812462) for the ``\emph{Macaulay2} Workshop at Wisconsin" on April 14-17, 2018, during which most of the work here was carried out. MB was supported by NSF award DMS-1502553; RL was supported by NSF award DMS-1502553; MP was supported by the NSF Graduate Research Fellowship (NSF award DGE-1313583); CR was supported by a Sloan fellowship and by NSF DMS-1600765; and SS was partially supported by a Sloan fellowship, NSF DMS-1500069, and NSF DMS-1651327.}
\begin{abstract}
We describe a \emph{Macaulay2} package for computing Schur complexes. This package expands on the \texttt{ChainComplexOperations} package by David Eisenbud. 
\end{abstract}

\section{Introduction}

Let $R$ be a commutative ring. The goal of this article is to describe the \emph{Macaulay2} package \texttt{SchurComplexes}, which computes the Schur complex $\S_{\l}(F)$ associated to a bounded complex $F$ of finitely generated free $R$-modules and a partition $\l$. 

Schur complexes are a simultaneous generalization of the symmetric and exterior power operations on complexes. The notion of a Schur complex was introduced by Nielsen in \cite{nielsen} in the characteristic 0 setting, and it was generalized to the characteristic free setting by Akin--Buchsbaum--Weyman in \cite{ABW}. 

The importance of such operations on complexes is illustrated by Walker's recent proof of the weak Buchsbaum--Eisenbud--Horrocks conjecture \cite[Theorem 2.4]{walker} in which exterior and symmetric squares of complexes play a crucial role. Walker's breakthrough work led to Eisenbud's implementation of the second exterior and symmetric power for complexes in the \emph{Macaulay2} package \texttt{ChainComplexOperations}. The work suggests that properties of Schur complexes should be further developed. Our goal here is to expand on this package by implementing the construction of an arbitrary Schur complex.

In Section \ref{background}, we provide some background on Schur complexes, following the detailed treatment in Weyman's book \cite[Section 2.4]{weyman}. In particular, we recall the ``straightening algorithm" of \cite{ABW} which expresses a $\Z/2$-graded tableau as a $\Z$-linear combination of so-called ``standard" $\Z/2$-graded tableaux; the implementation of this algorithm is the key component of the \texttt{SchurComplexes} package. Section \ref{examples} contains some examples of computations using \texttt{SchurComplexes}.

\section{Background on Schur complexes}
\label{background}
This section closely follows \cite[Section 2.4]{weyman}. Let $F = (F_0 \leftarrow F_1 \leftarrow \cdots \leftarrow F_d)$
be a bounded complex of finitely generated free $R$-modules. Denote by $\Feven$ (resp. $\Fodd$) the direct sum of the even (resp. odd) degree components of $F$, and choose bases $\{e_1, \dots, e_m\}$ and $\{f_1, \dots, f_n\}$ of $\Fodd$ and $\Feven$, respectively, which are unions of bases of the $F_i$.

\subsection{Exterior powers of complexes}

Fix a positive integer $r$, and let $T^r(F)$ denote the $r^{\th}$ tensor power of $F$. $T^r(F)$ may be equipped with a $\Sigma_r$-action in the following way: 
$$
\sigma \cdot (x_1 \otimes \cdots \otimes x_r) = \pm x_{\sigma^{-1}(1)} \otimes \cdots \otimes x_{\sigma^{-1}(r)},
$$
where the $x_i$ are homogeneous elements of $F$, and the sign is determined by declaring that transposing the elements $x_i$ and $x_j$ contributes the sign $(-1)^{|x_i||x_j|}$. By \cite[page 74]{weyman}, the $\Sigma_r$-action is compatible with the differential on $T^r(F)$. Let
$$
\varepsilon\colon T^r(F) \to T^r(F)
$$ 
denote the $R$-linear anti-symmetrization map
$$
(x_1 \otimes \cdots \otimes x_r) \mapsto \sum_{\sigma \in \Sigma_{r}} (-1)^{\on{sign}(\sigma)} \sigma \cdot (x_{1} \otimes \cdots \otimes x_r ).
$$

\begin{defn}
\label{exteriordefinition}
The $r^{\th}$ exterior power $\bigwedge^r F$ is defined to be the subcomplex $\varepsilon(T^r(F))$ of $T^r(F)$.
\end{defn}

\begin{rem}
\label{quotient}
Suppose $F$ is concentrated in even degrees, and let $I_r$ denote the $R$-submodule of $T^r(F)$ spanned by elements of the form
$$
x_1 \otimes \cdots \otimes x \otimes x \otimes \cdots \otimes x_r.
$$
By \cite[Section 1.2]{ABW}, $\varepsilon$ determines a split injection $T^r(F)/I_r \to T^r(F)$, and so Definition \ref{exteriordefinition} recovers the usual definition of an exterior power in this case. As noted in loc. cit., $\bigwedge^r F$ coincides with the antisymmetric tensors in $T^r(F)$ (i.e. those elements $v \in T^r(F)$ such that $\sigma \cdot v = (-1)^{\on{sign}(\sigma) }v$ for all $\sigma \in \Sigma_r$) when $\on{char}(R) \ne 2$. 
\end{rem}

There is a canonical $R$-linear embedding
$$
\iota \colon \bigoplus_i D_i(\Fodd) \otimes (T^{r-i}(\Feven)/I_{r-i}) \into T^r(F)
$$
whose image is precisely $\bigwedge^r F$. Here, $D_i(\Fodd) := \Sym_i((\Fodd)^*)^*$,
the $i^{\on{th}}$ divided power of $\Fodd$ (where $( - )^*$ denotes the $R$-linear dual). We now describe this embedding in detail. 

Let $I$ denote the two-sided ideal of $T(\Feven)$ generated by elements of the form $x \otimes x$, where $x \in F$. By Remark \ref{quotient}, the composition $T(\Feven) \into T(F) \xra{\varepsilon} T(F)$ factors through $T(\Feven)/I$ and induces an embedding
$$
\iota_\Lambda \colon T(\Feven)/I \into T(\Feven).
$$
By \cite{roby} Proposition IV.5, there is an embedding of $R$-algebras
$$
\iota_D\colon D(\Fodd) \to T(\Fodd)
$$
such that $\iota_D(x^{(j)}) =  \underbrace{x \otimes \cdots \otimes x}_{j \text{ copies}}$ for all $x \in \Fodd$, where the target is equipped with the shuffle product. The embedding $\iota$ is defined on each summand $D_i(\Fodd) \otimes (T^{r - i}(\Feven)/I_{r-i})$ by  
$$
\sum_{\sigma \in \Sigma_{i, r-i}} (-1)^{\on{sign}(\sigma)} \sigma \cdot (\iota_D \otimes \iota_\Lambda),
$$
where $\Sigma_{i, r-i} \subseteq \Sigma_r$ denotes the set of $(i, r-i)$ shuffles.

\begin{ex}
\label{koszul}
Let $x, y \in R$, and take $F$ to be the Koszul complex
$$
R \xla{\begin{pmatrix} x & y \end{pmatrix}} R^{\oplus 2} \xla{\begin{pmatrix} -y \\ x \end{pmatrix}} R
$$
on $x$ and $y$, lying in homological degrees $0,1$, and 2. Then $\bigwedge^2(F)$ is the complex
$$
	R^{\oplus 2} \xla{
	\begin{pmatrix}
		y&x&0&x\\
		0&y&x&-y\\
      	\end{pmatrix}
	} 
      R^{\oplus 4} \xla{
      \begin{pmatrix}
      		2 x& 0\\
      		{-y}&x\\
      		0&{-2 y}\\
      		{-y}&{-x}\\
      \end{pmatrix}
      } 
      R^{\oplus 2},
$$
lying in homological degrees $1,2$, and 3.
\end{ex}

The complex $\bigwedge F := \bigoplus_{r \ge 0} \bigwedge^r F$ is equipped with a product
$$\mu \colon \bigwedge^{r_1} F \otimes \bigwedge^{r_2} F \to \bigwedge^{r_1 + r_2} F$$
and a coproduct 
$$\Delta\colon \bigwedge^r F \to \bigoplus_{r_1 + r_2 = r}  \bigwedge^{r_1} F \otimes \bigwedge^{r_2} F.$$
For explicit formulas for $\mu$ and $\Delta$, we refer the reader to the proof of \cite[Proposition 2.4.1]{weyman}.

\subsection{Schur complexes}
Let $r$ be a positive integer, and let $\l = (\l_1, \dots, \l_s)$ be a partition of $r$, where $\l_i \ge \l_{i +1 }$. We will encode partitions with Young diagrams. For example, the partition $(3, 2, 2)$ of $7$ corresponds to the diagram
$$
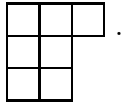
\ytableausetup
{mathmode, boxsize=1em}
\begin{ytableau}
 \text{ } &  &  \\
 &    \\
 &  \\
\end{ytableau} \text{ .}
$$
Let $c_1, \dots, c_t$ denote the lengths of the columns of $\l$.

\begin{defn}
\label{schurdefn}
The \emph{Schur complex} $\S_\l(F)$ is the quotient $(\bigwedge^{c_1} F \otimes \cdots \otimes \bigwedge^{c_t} F) / R$, where $R$ is the sum of submodules
$$
\bigwedge^{c_1} F \otimes \cdots \otimes \bigwedge^{c_{a - 1}} F \otimes R_{a, a+1} \otimes \bigwedge^{c_{a + 2}}F \otimes \cdots \otimes \bigwedge^{c_{t}}F.
$$
Here, $R_{a, a+1}$ is the submodule spanned by the images of the compositions
\begin{align*}
\Theta(a,u,v ; F)\colon \bigwedge^u F \otimes \bigwedge^{c_a - u + c_{a+1} - v} F \otimes \bigwedge^v F 
\xra{1 \otimes \Delta \otimes 1} & 
\bigwedge^u F \otimes \bigwedge^{c_a - u} F \otimes  \bigwedge^{c_{a+1} - v} F \otimes \bigwedge^v F  \\
\xra{\mu \otimes \mu}  \bigwedge^{c_a} & F \otimes  \bigwedge^{c_{a + 1}} F
\end{align*}
for $u + v < c_{a + 1}$.
\end{defn}

\begin{rem}
\label{discrepancy}
Our definition of the Schur complex differs from the one in \cite[Section 2.4]{weyman} in that the roles of the rows and columns are swapped. In other words, Weyman's definition of the Schur complex of $F$ with respect to $\l$ is recovered by applying Definition \ref{schurdefn} to the complex $F$ and the conjugate partition $\l^*$, i.e. the result of transposing the rows and columns of $\l$. 
\end{rem}

\begin{ex}
Of course, if $\l = (1, \dots, 1)$, $\S_\l(F) = \bigwedge^r F$. If $\l = (r)$, $\S_\l(F) = \Sym_r(F)$ (see \cite[Section 2.4]{weyman} for the definition of the symmetric power of a complex). 
\end{ex}

The basis of $F$ chosen above determines a basis $\mathcal{B}_i$  for each $\bigwedge^{c_i} F$: namely, the images under $\iota$ of elements of the form
$$
\{ 
e_1^{(j_1)} \cdots e_m^{(j_l)} \otimes f_{k_1} \wedge \cdots \wedge f_{k_{c_i - (j_1 + \cdots + j_l)}}
\}
$$
(from now on, we will tacitly identify these elements with their images under $\iota$). The set $\mathcal{S} := \{b_1 \otimes \cdots \otimes b_t \text{ : } b_i \in \mathcal{B}_i \}$ therefore gives an $R$-linear spanning set for $\S_\l(F)$. We will write elements of $\mathcal{S}$ as \emph{$\Z/2$-graded Young tableaux of shape $\l$}, i.e. functions
$$T\colon \{1, \dots, r \} \to \{-m, \dots, -1 \} \cup  \{1, \dots, n \} ,$$
where, as above, $m = \rk(\Fodd)$ and $n = \rk(\Feven)$. Here, divided power factors correspond to negative values, and exterior factors correspond to positive values. For instance, if $\lambda = (3,3,1)$ and $m = 2 = n$, the element
$(e_2^{(2)} \otimes f_1) \otimes (e_1 \otimes f_1) \otimes (f_1 \wedge f_2)$ in $\S_{(3,3,1)}(F)$ corresponds to the function

$$
\begin{tabular}{ c | c | c | c | c | c | c | c   }		
$l$  & 1 & 2 & 3 & 4 & 5 & 6 & 7  \\
   \hline
$T(l)$ &  -2 & -2 & 1 & -1 & 1 & 1 & 2
\end{tabular},
$$
which we express as the following Young tableau:

$$
T = \ytableausetup
{mathmode, boxsize=2em}
\begin{ytableau}
-2 & -1  & 1\\
-2 & 1  & 2 \\
1 \\
\end{ytableau} \text{ .}
$$
We will call such a tableau \emph{standard} if
\begin{itemize}
\item[(A)] the columns increase from top to bottom, with equality possible only for negative values, and
\item[(B)] the rows increase from left to right, with equality possible only for positive values.
\end{itemize}

\begin{rem}
  Our definition of a standard tableau is the transpose of Weyman's in \cite[Definition 1.1.12(c)]{weyman}   (cf. Remark~\ref{discrepancy} above).
\end{rem} 

For instance, the tableau $T$ above is standard. The tableaux
$$
\ytableausetup
{mathmode, boxsize=2em}
\begin{ytableau}
-2 & -1  & 1\\
-2 & -1  & 2 \\
-3 \\
\end{ytableau} \text{ ,}
\qquad
\ytableausetup
{mathmode, boxsize=2em}
\begin{ytableau}
-2 & -1  & 1\\
1 & 1  & 2 \\
1 \\
\end{ytableau} \text{ ,}
\qquad
\ytableausetup
{mathmode, boxsize=2em}
\begin{ytableau}
-2 & -1  & -1\\
-2 & 1  & 2 \\
1 \\
\end{ytableau}
$$
are non-standard. 
\begin{prop}[{\cite[Proposition 2.4.2]{weyman}}]
The standard $\Z/2$-graded tableaux of shape $\l$ form an $R$-linear basis of $\S_\l(F)$.
\end{prop}

We compute the differentials in $\S_\l(F)$ with respect to this basis in the \texttt{SchurComplexes} package. It is therefore essential for us to implement an algorithm for writing the image of a standard tableau under the differential in $\S_\l(F)$ as a linear combination of standard tableaux. The proof of \cite[Proposition 2.4.2]{weyman} explains such an algorithm: the ``straightening algorithm" of \cite{ABW}.  We now discuss this algorithm in detail.

\subsection{The straightening algorithm}

As in the previous subsection, $\l = (\l_1, \dots, \l_s)$ denotes a partition of a positive integer $r$ with column lengths $c_1, \dots, c_t$. The straightening algorithm is a process for writing a $\Z/2$-graded tableau in the spanning set $\mathcal{S}$ of $\S_\l(F)$ described above as a $\Z$-linear combination of standard tableaux. Here is how it works. 

Let $T$ be a tableau in $\mathcal{S}$. We make the following observations:
\begin{itemize}
\item If a column of $T$ contains a repeated positive entry, $T = 0$ in $\S_\l(F)$.
\item Since the divided power (resp. exterior) algebra is commutative (resp. skew commutative), rearranging the columns in $T$ so that it satisfies (A) only changes the element of $\S_\l(F)$ represented by $T$ up to a sign.
\end{itemize}

With these facts in mind, we recall the straightening algorithm:

\begin{itemize}
\item{{\bf Input}:} A tableau $T \in \mathcal{S}$.
\vskip.1in
\item{{\bf Step 1:}} Denote by $T'$ the result of rearranging the columns of $T$ so that they satisfy (A), and let $\sigma \in \{\pm 1\}$ denote the resulting sign.
\vskip.1in
\item{{\bf Step 2:}} If $T'$ satisfies (B), output $\sigma T'$. Otherwise, choose the topmost row, say the $w^{\on{th}}$ row, with a ``violation" of (B). Let $T'(i,j)$ denote the entry with horizontal coordinate $i$ and vertical coordinate $j$ in the Young diagram, counting from the top-left corner. So $T'(i,j)$ is the entry in the $i^{\on{th}}$ column and the $j^{\on{th}}$ row. Choose the smallest index $a$ such that either $T'(a, w) > T'(a+1 ,w) $ or $T'(a,w) = T'(a+1, w) < 0$. Then, choose the smallest index $w'$ such that $T'(a+1, w) < T'( a+1, w' + 1)$; if no such index exists, set $w' = c_{a+1}$.   
\vskip.1in
\item{{\bf Step 3:}} Set $u := w -1$ and $v := c_{a + 1} - w'$. Define
\begin{enumerate}
\item $V_1 \in \bigwedge^u F$ to be the element corresponding to the first $u$ entries in the $a^{\on{th}}$ column of $T'$,
\item $V_2 \in \bigwedge^{c_a - u + c_{a + 1} - v} F$ to be the element corresponding to the last $c_{a} - u$ entries in the $a^{\on{th}}$ column of $T'$ followed by the first $c_{a+1} - v = w'$ entries of the $(a+1)^{\on{st}}$ column of $T'$,
\item $V_3 \in \bigwedge^v F$ to be the element corresponding to the last $v$ entries in the $(a+1)^{\on{st}}$ column of $T'$.
\end{enumerate}
For $k \in \{1, \dots, a-1, a+2, \dots, s\}$, define $U_k$ to be the element of $\bigwedge^{c_k} F$ which corresponds to the $k^{\on{th}}$ column of $T'$. Recall that
$$
L:=(1 \otimes \Theta(a,u,v; F) \otimes 1)
(U_1 \otimes \cdots \otimes U_{a - 1} \otimes V_1 \otimes V_2 \otimes V_3 \otimes U_{a+2} \otimes \cdots \otimes U_{s})
$$
is 0 in $\S_\l(F)$, where $\Theta(a,u,v; F)$ is as in Definition \ref{schurdefn}. $L$ is a $\Z$-linear combination of tableaux in $\mathcal{S}$ which contains $T'$ with coefficient 1. If each tableau in the sum $T' -L$ is standard, output $\sigma(T' - L)$. Otherwise, repeat this algorithm on each tableaux in $\sigma(T' - L)$, keeping track of the coefficients.
\end{itemize}

The key observation is that each tableau appearing in the linear combination $T'-L$ from Step 3 is strictly ``smaller" than $T'$, in the sense described in \cite[Section 1.1]{weyman}, and so the algorithm does indeed terminate. 

\begin{ex}
\label{straighteningexample}
Let's apply the straightening algorithm to the tableau 
$$
T = \ytableausetup
{mathmode, boxsize=2em}
\begin{ytableau}
-3 & 2  & -1\\
-2 & 1  & 3\\
-2 & 3 \\
\end{ytableau} \text{ .}
$$

\begin{itemize}
\item{{\bf Step 1:}} The middle column needs to be rearranged. Since $f_2f_1f_3 = -f_1f_2f_3$ in $\bigwedge^3 \Feven$, we have 
$$
T' = \ytableausetup
{mathmode, boxsize=2em}
\begin{ytableau}
-3 & 1  & -1\\
-2 & 2  &  3\\
-2 & 3 \\
\end{ytableau} 
$$
and $\sigma = -1$.
\vskip.1in
\item{{\bf Step 2:}} $T'$ is not standard. Here, $w = 1$, $a = 2$, and $w' = 1$. 
\vskip.1in
\item{{\bf Step 3:}} Here, $u = 0$ and $v = 1$. We have $U_1 = e_3e_2^{(2)}$, $V_2 = e_1 \otimes f_1 \wedge f_2 \wedge f_3$, and $V_3 = f_3$ (since $u = 0$, $V_1$ plays no role), so 
$$
L = (1 \otimes \Theta(2, 0, 1; F))(e_3e_2^{(2)} \otimes (e_1 \otimes f_1 \wedge f_2 \wedge f_3) \otimes f_3).
$$
In this case, $\Theta(2, 0, 1)$ is the composition
$$
(\bigwedge^4 F) \otimes F \xra{\Delta \otimes 1}   (\bigwedge^3 F \otimes  F)  \otimes F 
\xra{\id \otimes \mu}  \bigwedge^3 F \otimes  \bigwedge^2 F.
$$
By the proof of \cite[Proposition 2.4.1(a)]{weyman}, the relevant component of the coproduct 
$$
\Delta\colon \bigwedge^4 F \to \bigwedge^3 F \otimes F
$$ 
in $\Theta(2, 0, 1; F)$ is:
\begin{align*}
D_1(\Fodd) \otimes \bigwedge^3 \Feven \xra{\Delta_D \otimes \Delta_\Lambda} & 
\bigoplus_{i = 0}^1 D_i(\Fodd) \otimes  D_{1 - i}(\Fodd) \otimes \bigwedge^{3 - i} \Feven \otimes \bigwedge^i \Feven  \\
\xra{\tau} &
\bigoplus_{i = 0}^1 D_i(\Fodd) \otimes   \bigwedge^{3 - i}  \Feven \otimes D_{1 - i}(\Fodd) \otimes \bigwedge^i \Feven,
\end{align*}
where $\Delta_D$ and $\Delta_\Lambda$ are the coproducts on the ordinary divided power and exterior algebra, and $\tau$ transposes the middle factors and multiplies by the sign $(-1)^{(1-i)(3 - i)}$.

Note: one might think that, since the elements of $\bigwedge^{3 - i} \Feven$ have even degree, transposing the middle two factors should not introduce a sign. But recall that we are considering $D_{1 - i}(\Fodd) \otimes \bigwedge^{3 - i} \Feven$ as a submodule of $T^4(F)$ via the embedding $\iota$, which shuffles together the elements of $D_{1 - i}(\Fodd)$ and $\bigwedge^{3 - i} \Feven$: this is why it is necessary to multiply by $(-1)^{(1-i)(3 - i)}$. The general rule here is: when one transposes the factors of $D_{s}(\Fodd) \otimes \bigwedge^{t} \Feven$, one must introduce the sign $(-1)^{st} = \on{sign}((1 \text{ } 2 \cdots s + t )^t$).

Applying the formula for $\Delta_\Lambda$ in \cite[Section 1.1, page 3]{weyman}, one gets
\begin{align*}
L = -(e_3e_2^{(2)}) \otimes (f_1 \wedge f_2 \wedge f_3) \otimes (e_1 \otimes f_3) & - (e_3e_2^{(2)}) \otimes (e_1 \otimes f_1 \wedge f_3) \otimes (f_2 \wedge f_3)  \\
&+ (e_3e_2^{(2)}) \otimes (e_1 \otimes f_2 \wedge f_3) \otimes (f_1 \wedge f_3),\\
\end{align*}
and therefore
$$
T =  
\ytableausetup
{mathmode, boxsize=2em}
\begin{ytableau}
-3 & -1  & 2 \\
-2 &  1  &  3\\
-2 & 3 \\
\end{ytableau} 
-\ytableausetup
{mathmode, boxsize=2em}
\begin{ytableau}
-3 & -1  & 1 \\
-2 &  2  &  3\\
-2 & 3 \\
\end{ytableau} \text{ .}
$$
Both of these tableaux are standard, so we're done. 
\end{itemize}
\end{ex}

\section{Examples of computations using the package \texttt{SchurComplexes}}
\label{examples}
The \texttt{SchurComplexes} package has two main functions:
\begin{itemize}
\item {\tt straightenTableau}, which applies the straightening algorithm to a tableau, and
\item {\tt schurComplex}, which computes the Schur complex of a bounded complex of finitely generated free modules.
\end{itemize}
\subsection{Using the function {\tt straightenTableau} }
Let's apply {\tt straightenTableau} to the tableau in Example \ref{straighteningexample}. First, we load the package:

\begin{verbatim}
i1 : loadPackage "SchurComplexes.m2"

\end{verbatim}

We encode the tableau $T$ in a hash table:

\begin{verbatim}
i2 : T = new HashTable from {(1,1) => -3, (1,2) => -2, (1,3) => -2, 
(2,1) => 2, (2,2) => 1, (2,3) => 3, (3,1) => -1, (3,2) => 3}

\end{verbatim}



We encode the partition $(3,3,2)$ in a list:

\begin{verbatim}
i3 : lambda = {3,3,2}

\end{verbatim}



Now, we apply {\tt straightenTableau} to the pair {\tt(T, lambda)}:

\begin{verbatim}
i4 : straightenTableau(T, lambda)

o4 = HashTable{HashTable{(1, 1) => -3} => 1 }
                         (1, 2) => -2
                         (1, 3) => -2
                         (2, 1) => -1
                         (2, 2) =>  1
                         (2, 3) =>  3
                         (3, 1) =>  2
                         (3, 2) =>  3
               HashTable{(1, 1) => -3} => -1
                         (1, 2) => -2
                         (1, 3) => -2
                         (2, 1) => -1
                         (2, 2) =>  2
                         (2, 3) =>  3
                         (3, 1) =>  1
                         (3, 2) =>  3

\end{verbatim}



The output is a hash table which assigns a coefficient to each standard tableau in the linear combination comprising the straightening of $T$. Notice that the output agrees with the calculation in Example \ref{straighteningexample}.

\subsection{Using the function {\tt schurComplex}}
Let $R=\mathbb{Q}[x_{i,j}]$, where $1 \le i \le 2$ and $1 \le j \le 4$, and let
$$F\colon R^{\oplus 4} \xra{(x_{i,j})} R^{\oplus 2}$$
denote the generic $2 \times 4$ matrix, considered as a complex concentrated in degrees 1 and 0. By 
\cite[Exercise 6.34(d)]{weyman}, $\S_{(3)}(F) = \Sym_3(F)$ has nonzero homology only in degree 0. We now use the {\tt schurComplex} function to compute $\S_{(3)}(F)$ and verify this fact.

We first load the package, fix our ground ring $R$, and define our complex $F$: 
\begin{verbatim}
i1 : loadPackage "SchurComplexes.m2"
i2 : R = QQ[x11,x21,x12,x22,x13,x23,x14,x24];
i3 : M = genericMatrix(R,x11,2,4);
i4 : F = new ChainComplex; 
i5 : F.ring = R; F#0 = target M; F#1 = source M; F.dd#1 = M;

\end{verbatim}

The function {\tt schurComplex} takes as input a {\tt ChainComplex} and a {\tt List} which encodes the partition. Let's define our partition and compute $\S_{(3)}(F)$:

\begin{verbatim}
i6 : lambda = {3};
i7 : S = schurComplex(lambda,F)

      4      12      12      4
o7 = R  <-- R   <-- R   <-- R
                              
     0      1       2       3

\end{verbatim}

Finally, let's check that $\S_{(3)}(F)$ has trivial homology in degrees greater than 0:
\begin{verbatim}
i8 : apply((length S)+1,i->reduceHilbert hilbertSeries HH_i(S))

          4     0  0  0
o8 = {--------, -, -, -}
             5  1  1  1
      (1 - T)

\end{verbatim}

\vskip\baselineskip
\noindent {\bf Acknowledgements.} We thank Mike Stillman for helpful comments while the  \texttt{SchurComplexes} package was being written, and we thank the University of Wisconsin-Madison Mathematics Department for hosting the NSF-funded ``\emph{Macaulay2} Workshop at Wisconsin" on April 14-17, 2018, during which most of this package was written. We also thank the referees for their helpful comments on a previous version of this article.
\bibliographystyle{amsalpha}
\bibliography{Bibliography}
\Addresses
\end{document}